\documentclass[11pt,a4paper]{article}
\usepackage{amsmath}
\usepackage{amsthm}
\usepackage{amssymb}
\usepackage{graphics}
\usepackage{colortab} 
\usepackage{enumitem}

\theoremstyle{plain}
  \newtheorem{theorem}{\bf Theorem}[section]
  
  \newtheorem{lemma}[theorem]{\bf Lemma}

\theoremstyle{remark}
  \newtheorem{remark}[theorem]{\bf Remark}

\begin{document}
\title{Uniqueness for the inverse fixed angle scattering problem.}
\author{J. A. Barcel\'{o}\thanks{Departamento de Matem\'atica e Inform\'atica aplicadas a las Ingenier\'ias Civil y Naval, Universidad Polit\'ecnica de Madrid. E-mails: juanantonio.barcelo@upm.es, carlos.castro@upm.es, maricruz.vilela@upm.es.} \and C. Castro\footnotemark[1] \and  T. Luque\thanks{Departamento de An\'alisis Matem\'atico, Facultad de Matem\'aticas, Universidad Complutense de Madrid. E-mail: t.luque@ucm.es.} \and C. J. Mero\~no\footnotemark[1]  \and A. Ruiz\thanks{Departamento de Matem\'aticas, Facultad de Ciencias, Universidad Autónoma de Madrid. E.mail: alberto.ruiz@uam.es} \and M. C. Vilela\footnotemark[1]}
\providecommand{\keywords}[1]{\textbf{\text{Key words.}} #1}
\maketitle

\abstract{We present a uniqueness result in dimensions $2$ and $3$ for the inverse fixed angle scattering problem associated to the Schr\"o\-din\-ger operator $-\Delta+q$, where $q$ is a small real valued potential with compact support in the Sobolev space $W^{\beta,2}$ with $\beta>0.$ 
	This result improves the known result, due to Stefanov \cite{St}, in the sense that almost no regularity is required for the potential.
	The uniqueness result still holds in dimension $4$, but for more regular potentials in $W^{\beta,2}$ with $\beta>2/3$. 
	
	The proof is a consequence of the reconstruction method presented in our previous work \cite{BCLV}.}

\vspace{0.5cm}

\section{Introduction}
We consider the scattering problem appearing in quantum mechanics of finding 
$u$ the solution of
\begin{equation}\label{u_problem}
\left\{  \begin{array}{lll}
(-\Delta +q(x)-k^2)u(x)=0, \hspace{1cm} x \in \mathbb{R}^d, \\
u=u_i+u_s, \\ \partial_ru_s-iku_s=o\left( r^{-(d-1)/2}  \right), \hspace{0.5cm} r=|x|\longrightarrow \infty,

\end{array}      \right.
\end{equation}
where $q\equiv q(x)$ is a real valued potential with compact support,
$u_i\equiv u_i(x, \theta, k)=e^{ik \theta \cdot x}$ is the incident wave, with wave number $k>0$ and direction of propagation $\theta \in \mathbb{S}^{d-1}$. 

The solution $u$ of problem \eqref{u_problem} is known as the outgoing scattering solution, and $u_s$, which is the perturbation of the incident wave due to the potential, as the scattered wave.

It is well known that for appropriate potentials $q$, $u_s$ satisfies
\begin{equation}\label{u_asintotica}
u_s(x, \theta, k)=C_d k^{(d-3)/2}\frac{e^{ik |x|}}{|x|^{(d-1)/2}}u_\infty (\omega, \theta,k)+o\left( |x|^{-(d-1)/2}  \right),
\quad
|x|\rightarrow\infty,
\end{equation}
where $\omega= x/|x|$ and 
\begin{equation}\label{u_amplitud}
u_\infty (\omega, \theta ,k)=\int_{\mathbb{R}^d}e^{-ik \omega \cdot y}q(y)u(y,\theta, k)dy.
\end{equation}

The function $u_\infty\equiv u_\infty (\omega,\theta,k)$, known as the scattering  amplitude  or far-field pattern of $u$, depends on the wave number $k>0$, the incident angle $\theta\in\mathbb{S}^{d-1}$  and the reflecting angle $\omega \in\mathbb{S}^{d-1}$. The central problem in inverse scattering for the Schr\"{o}dinger
equation is to recover the unknown potential $q$ from the scattering amplitude $u_\infty$. This is an overdetermined problem, since  
$u_\infty$ depends on $2d-1$ variables  while $q$ only depends on $d$.

Therefore, it is natural to recover the potential from partial knowledge of the scattering amplitude. From the applications point of view, the most widely studied cases are the fixed angle scattering problem, the fixed energy scattering problem and the backscattering one. An important question that arise is if this partial knowledge of the scattering amplitude determines the potential. Here we are interested in this problem when considering fixed angle scattering data, i.e. given $\theta_0\in \mathbb{S}^{d-1}$, we assume known $u_\infty\equiv u_\infty (\omega,\pm \theta_0,k)$ for $\omega\in\mathbb{S}^{d-1}$ and $k>0$.  More precisely, we are interested in the following uniqueness question: Are there two different potentials $q_1$ and $q_2$ with the same fixed angle scattering data? Note that in this case both the potential and the scattering data depend on the same number of parameters $d$.

Mathematical analysis of reconstruction algorithms for fixed angle scattering data was first considered by Prosser \cite{Pr}, who obtained a potential in series form that fits the fixed angle scattering data in dimension three, under some restrictions on the size of the potentials in a Friedrichs type norm. However, the uniqueness result remained open. The first partial uniqueness result is due to Bayliss,  Li and  Morawetz \cite{BLM}, that proved uniqueness when the fixed angle scattering amplitude vanishes. This result was extended by Ramm \cite{Ra} to scattering data coming from compactly supported potentials $q$ that are  constant along the lines defined by the direction $\theta=\theta_0$.

Later on, Stefanov \cite{St} proved the first uniqueness result for small complex valued potentials in $W^{4, \infty}$ with compact support in $\mathbb{R}^3$. In the particular case of real potentials, he also showed that it is enough to consider scattering data for a unique incident angle $\theta_0$ (see Remark \ref{RE}). 

For real compact support potentials in $\mathbb{R}^3$, Ramm \cite{Ra2} stated a global uniqueness result in $W_0^{4,2}$ without requiring any smallness condition. However, there is a gap in his proof. The mistake appears in an essential part of the argument and it is not clear how to avoid it. We refer to the Appendix below for more details. Therefore the global uniqueness for the fixed angle scattering problem is still open. Very recently, a global uniqueness result using the time dependent approach to scattering  has been given   in \cite{RS19} for smooth potentials with compact support. 

Our uniqueness result, given in the following theorem, improves the  above mentioned results since much less regularity is required on the potentials. However, as in \cite{RS19}, we consider a larger set of scattering data since we use the scattering amplitude coming from both $\theta_0$ and $-\theta_0$.

\begin{theorem}	\label{u_teorema}
	Let $\beta >0$ if $d=2$ or $3$ and $\beta>2/3$ if $d=4$. For $i=1,2$ let $q_i$ be  a real valued potential in $W^{\beta , 2}(\mathbb{R}^d)$ with compact support and scattering amplitude given by $u^i_\infty (\omega, \theta , k)$. Then, there exists a constant $B$ small enough (see Remark \ref{u_nota_constante}), depending on $\beta, d$ and the supports of $q_1$ and $q_2$, such that if
	\begin{equation}\label{u_condicion_1_teorema}
	\|q_i\|_{W^{\beta,2}} < B, \;\;\; i=1,2,
	\end{equation}
	and 
	\begin{equation}\label{u_condicion_2_teorema}
	u^1_\infty(\omega, \pm\theta_0, k)=u^2_\infty(\omega, \pm\theta_0, k), \;\;\; \textrm{ for all } \; \omega \in \mathbb{S}^{d-1} \; \textrm{ and } k>0,
	\end{equation} 
	then $q_1=q_2$.
\end{theorem}

In general, $W^{\alpha , p}(\mathbb{R}^d)$ denotes the Sobolev space
\begin{equation*}
\label{sobolev}
W^{\alpha,p}(\mathbb{R}^d)=\{f\in\mathcal{S}^{\prime}(\mathbb{R}^{d}):\Lambda^{\alpha}f\in L^{p}(\mathbb{R}^d)\},\qquad \alpha\in\mathbb{R},1\le p\le\infty,
\end{equation*}
where
$$
\Lambda^{\alpha}
=
(1-\Delta)^{\alpha/2}
=
\mathcal{F}^{-1}\left\langle \xi\right\rangle^{\alpha}\mathcal{F},
$$
$\mathcal{F}$ denotes the Fourier transform and
$\left\langle\xi\right\rangle =(1+\left\vert \xi\right\vert ^{2})^{1/2}.$
\section{Proof of Theorem \ref{u_teorema}}
The proof of Theorem \ref{u_teorema} is a direct consequence of Theorem 1.1. of \cite{BCLV}, where authors construct a sequence of approximations using the fixed angle scattering data and prove that it converges to the potential in certain Sobolev space and for certain potentials.

For completeness, we include here the statement of Theorem 1.1. of \cite{BCLV}. We first observe that the scattered wave $u_s$ in \eqref{u_problem} for $\theta=\theta_0$ satisfies the Lippmann-Schwinger integral equation given by
\begin{equation}\label{u_L-S}
u_s(x, \theta_0,k)=R_k\left( qe^{ik \theta_0 \cdot (\cdot)}  \right)(x)+R_{k}\left( qu_s(\cdot, \theta_0,k )     \right)(x), \hspace{0.5cm} x \in \mathbb{R}^{d},
\end{equation}
where $R_k$ denotes the outgoing resolvent of the Laplace operator, which in terms of the Fourier transform, is defined by
\begin{equation}
\label{resolvente}
\widehat{R_k(f)}(\xi)=\frac{\widehat{f}(\xi)}{-|\xi|^2+k^2+i0}.
\end{equation}
Under the conditions of Theorem \ref{u_teorema}, the existence and uniqueness of $u_s$ is guaranteed by Theorem 2.11 of \cite{BCLV}.

Moreover, since $u(x, \theta_0,k)=e^{ik\theta_0\cdot x}+u_s(x, \theta_0,k)$, from \eqref{u_amplitud} we obtain
\begin{align}\label{u_1aproximacion}
u_\infty(\omega, \theta_0, k)=\widehat{q}(k(\omega-\theta_0))+\int_{\mathbb{R}^d}
e^{-ik\omega \cdot y}q(y)u_s(y,\theta_0,k) dy.
\end{align}
This provides a first approximation to the unknown potential $q$ called the Born approximation.
Formally, if we remove the last term in \eqref{u_1aproximacion}, the right hand side is the Fourier transform of $q$ along the so-called Ewald spheres. More precisely, given $\theta_0$ fixed, we have, up to a zero measure set,
$$\mathbb{R}^d=H_{\theta_0} \cup  H_{-\theta_0}=\left\{  \xi \in \mathbb{R}^d: \xi \cdot \theta_0 <0 \right\} \cup 
\left\{  \xi \in \mathbb{R}^d: \xi \cdot \theta_0 >0 \right\}.$$
For $\xi \in H_{\pm \theta_0}$ there exist unique $\omega (\xi) \in \mathbb{S}^{d-1}$ and $k(\xi)>0$ such that 
$$\xi=k(\xi)(\omega (\xi)-\theta_0(\xi)), \hspace{0.3cm} \textrm{ with } \; \theta_0(\xi)=\left\{  \begin{array}{ll}
\theta_0, \hspace{0.6cm} \xi \in H_{\theta_0} \\  -\theta_0, \hspace{0.3cm} \xi \in H_{-\theta_0} 
\end{array}    \right. .$$
Actually, if $\xi \in H_{\theta_0}$ we have that
\begin{equation}\label{u_lola}
\omega(\xi)=-\frac{2 \xi \cdot \theta_0}{|\xi|^2}\xi+\theta_{0},  \;\; \; \textrm{ and } \; \; k(\xi)=-\frac{|\xi|^2}{2 \xi \cdot \theta_0}.
\end{equation}

The Born approximation for fixed angle scattering data $\theta_0$ of a potential $q$ is defined by
\begin{equation}\label{u_aproximacion_born}
\widehat{q_{\theta_0}}(\xi)=u_\infty( w(\xi), \theta_0(\xi), k(\xi)  ).
\end{equation}
\begin{remark}\label{RE} Definition \eqref{u_aproximacion_born} is equivalent to the one used by Stefanov in \cite{St} when $\xi \in H_{\theta_0}$, but it differs when $\xi \in H_{-\theta_0}$. More precisely, instead of considering scattering data for $-\theta_0$ and $k>0$,  the scattering amplitude in \cite{St} is extended to $k<0$ by $u_{\infty}(w,\theta_0,k) = \overline{u_{\infty}(w,\theta_0,-k)}$. In this way, only scattering data coming from $\theta_0$ are required. 
	
\end{remark}

Now we insert iteratively the Lippmann-Schwinger integral equation (\ref{u_L-S}) into (\ref{u_1aproximacion}) to obtain for $\xi=k(\xi)(\omega (\xi)-\theta_0(\xi))$ and   $m\geq 1$
\begin{equation}\label{u_serie_born}
\widehat{q_{\theta_0}}(\xi)=\widehat{q}(\xi)+\sum_{j=1}^m \widehat{\mathcal{Q}_j(q)}(\xi)+\widehat{q_m^r}(\xi),
\end{equation}
where
\begin{equation}\label{u_serie_born_1}
\widehat{\mathcal{Q}_j(q)}(\xi)=\int_{\mathbb{R}^d}e^{-ik(\xi)\omega(\xi) \cdot y} (qR_{k(\xi)})^j \left( qe^{ik(\xi) \theta_0(\xi) \cdot (\cdot)}   \right)(y) dy,
\end{equation}
and 
$$\widehat{q_m^r}(\xi)=\int_{\mathbb{R}^d}e^{-ik(\xi)\omega(\xi) \cdot y} (qR_{k(\xi)})^m \left( qu_s\left(   \cdot, \theta_0(\xi),k(\xi) \right)  \right)(y) dy.$$

In order to construct our sequence of approximations of $q$, we fix $m$ and consider the operator
\begin{equation}\label{u_1_operador}
\mathcal{T}_m(f)=\phi q_{\theta_0}-\phi \sum_{j=1}^m \mathcal{Q}_j(f)(\xi),
\end{equation}
where $\phi \in \mathcal{C}^\infty $ is a cut-off function with compact support satisfying
$$\phi (x)=1, \; \textrm{ if } |x|<R \; \textrm{ and } \phi (x)=0, \; \textrm{ if } |x|>2R,$$
with $R>0$ such that the support of q is contained in $B(0,R)$.

We define the sequence $\left\{ q_{m,\ell}  \right\}_{\ell \in \mathbb{N}}$  recursively by
\begin{equation}\label{u_1_sucesion}
\left\{   \begin{array}{ll}
q_{m,1}=0, \\ q_{m,\ell +1}=\mathcal{T}_m(q_{m,\ell}), \hspace{0.4cm} \ell \geq 1.
\end{array}   \right.
\end{equation}
\begin{remark}
	\label{dependencia}
	Observe that $q_{m,2}$ is nothing but $\phi q_{\theta_0}$. 
	On the other hand, form (\ref{u_1_operador}), the sequence $\left\{ q_{m,\ell}  \right\}_{\ell \in \mathbb{N}}$ only depends  on $u_\infty(\omega, \pm\theta_0, k)$ and $R$.
\end{remark}

The following theorem states that the iterated limit, first in $\ell$ and then in $m$, of the sequence $\left\{ q_{m,\ell}  \right\}_{m,\ell \in \mathbb{N}}$ is equal to $q$, under certain hypotheses. 

\begin{theorem}\label{u_teorema_recuperacion}
	(Theorem 1.1 of \cite{BCLV}) For $d \geq 2$ and $\alpha $ satisfying 
	\begin{equation} \label{u_alba}
	0<\alpha \leq 1, \hspace{0.5cm} \textrm{ and} \hspace{0.5cm} \frac{d}{2}-\frac{d}{d-1}< \alpha <\frac{d}{2},
	\end{equation} 
	leq $q \in W^{\alpha,2}(\mathbb{R}^d)$ be a real function with compact support in $B(0,R)$ and such that 
	\begin{equation}\label{alba}
	\|q\|_{W^{\alpha, 2}}<A,
	\end{equation}
	for an appropriate constant $A:=A(d, \alpha , R)>0$ small enough. For each $m \in \mathbb{N}$, let $\left\{ q_{m,\ell}  \right\}_{\ell \in \mathbb{N}}$ be the sequence defined by (\ref{u_1_sucesion}). Them, there exists $q_m \in W^{\alpha,2}(\mathbb{R}^d)$ satisfying
	\begin{equation}\label{u_2_sucesion}
	q_m=\lim_{\ell \longrightarrow \infty}q_{m,\ell} \hspace{0.8cm} \textrm{in } W^{\alpha,2}(\mathbb{R}^d).
	\end{equation}
	Moreover, the sequence $\left\{ q_m\right\}_{m \in \mathbb{N}}$ satisfies
	\begin{equation}\label{u_3_sucesion}
	\lim_{m \longrightarrow \infty}q_{m}=q \hspace{0.8cm} \textrm{in } W^{\alpha,2}(\mathbb{R}^d).
	\end{equation}
\end{theorem}

\begin{remark}
	The proof of Theorem \ref{u_teorema_recuperacion} depends on very accurate estimates for $R_k$, the resolvent of the Laplace operator, which
	can be found in \cite{Ru1}. The smallness condition \eqref{alba} is written in terms of the constants appearing in these estimates (see Remark 1.2 of \cite{BCLV}). On the other hand, conditions in \eqref{u_alba} imply that $d<5$ and, therefore, Theorem \ref{u_teorema_recuperacion} is only valid for $2 \leq d \leq 4$. 
\end{remark}

\textit{Proof of Theorem \ref{u_teorema}}.
The proof is a straightforward consequence of Theorem \ref{u_teorema_recuperacion}. Precisely, since any $q$ under the conditions of Theorem \ref{u_teorema} satisfies the assumptions of Theorem \ref{u_teorema_recuperacion}, for $i=1,2$ we can construct the approximation sequence $\left\{ q^i_{m,\ell}  \right\}_{m,\ell \in \mathbb{N}}$ for $q_i$, which only depends on $u_\infty^i(x, \pm \theta_0,k)$ and the support of $q_i$. 
Therefore, by \eqref{u_condicion_2_teorema}, $q_{m,\ell}^1=q_{m,\ell}^2$. From (\ref{u_2_sucesion})-(\ref{u_3_sucesion}) we easily obtain $q_1=q_2$.
\hfill
$\Box$

\begin{remark} \label{u_nota_constante}
	Notice that in the particular cases $\beta <1$ and $d=2$ or $\beta \leq 1$ and $d=3$ or $d=4$, the constant $B$ given in \eqref{u_condicion_1_teorema} equals the constant $A$ given in \eqref{alba}, with $\alpha=\beta$. In other case, we can chose any $\alpha$ under the conditions of Theorem \ref{u_teorema} such that $0< \alpha < \beta$, and thus
	$$\|q_i \|_{W^{\alpha , 2}} \leq C(d, \alpha, \beta,R)  \|q_i \|_{W^{\beta ,2}}<A, \quad i=1,2,$$
	if we take
	$$B=\frac{A}{C(d, \alpha, \beta,R)}.$$
\end{remark}

\section{Appendix}
In this section we detail a gap in the proof of the global uniqueness for the fixed angle scattering problem given in \cite{Ra2}.  Essentially, the main result in \cite[Theorem 1.1]{Ra2} follows from putting together key estimates (25) and (26)  in that paper. 

Estimate (26) (which coincides with estimate \eqref{(26)} below) is  proved in  \cite[Lemma 3.1]{Ra2} that we now state.

\begin{lemma} Assume that $f$ is a real  function compactly supported in $B(0,R)$, and  such that $f\in W^{\beta, 2}(B(0,R))$ with $\beta>3$. Let $\kappa>0$ fixed and $\eta>0$, then
	\[\limsup_{\eta\to\infty}\max_{\zeta\in\mathbb{S}^2}|\hat f((\kappa+i\eta)\zeta)|=\infty.\]
	Also for any $\kappa>0$, there is an $\eta=\eta(\kappa)$ such that
	\begin{equation}\label{(26)}
	\max_{\zeta\in\mathbb{S}^2}|\hat f((\kappa+i\eta)\zeta)|\geq\max_{\xi\in\mathbb{R}^3}|\hat f(\xi)|,
	\end{equation}
	and 
	\[\eta(\kappa)=R^{-1}\ln\kappa+O(1),\quad \kappa\to\infty.\]
\end{lemma}


Our main objection is that \eqref{(26)} cannot be true for a  $C^\infty _0$ function and a sequence of the kind $\eta(\kappa)= O(\ln \kappa)$. In fact, if  $f \in C^\infty_0(B(0,R))$ by Paley-Wiener (see \cite[pp.21]{Hormander}) we have that 

\[ |\hat f((\kappa +i\eta)\zeta)|\leq C(\gamma) |\kappa+i \eta|^{-\gamma} e^{R|\eta|} \]
holds for  every $\gamma >0$. Then it is clear that  the condition  $\eta(\kappa)= O(\ln \kappa)$ as $\kappa\to \infty$, implies  that the left hand side in \eqref{(26)} tends to 0 as $\kappa\to \infty$, contradicting the inequality.

Since the logarithmic  behavior of $\eta(\kappa)$ is  a condition to prove estimate (25)  (see estimates (45)-(49) in \cite{Ra2}), then both estimates (25) and (26) cannot be used simultaneously together to prove  \cite[Theorem 1.1]{Ra2}.

\paragraph{Acknowledgements.} The first, the second, the fourth and the sixth author were supported by Spanish Grant MTM2017-85934-C3-3-P, the third by Spanish Grant MTM2017-82160-C2-1-P, and the fifth by Spanish Grant MTM2017-85934-C3-2-P.

\end{document}